\input amstex

\documentstyle{amsppt}
\magnification \magstep 1

\def\CC{\Bbb C}

\def\QQ{\Bbb Q}

\def\ot{\otimes}

\def\L{\Lambda}

\def\OO{\Cal{O}}

\def\ox{\omega_X}

\def\FF{\Cal{F}}

\def\UU{\Cal{U}}
\def\VV{\Cal{V}}

\def\LL{\Cal{L}}

\def\Pic{{\text{\rm Pic}}}
\def\PX{{\text {\rm Pic}}^0(X)}

\def\L1{\LL^{(1)}}

\def\WED#1#2#3{#1_{#2} \wedge \ldots \wedge #1_{#3}}
\def\DI#1{\text{\rm dim}#1}
\def\alb{\text{\rm alb}}
\def\Alb{\text{\rm Alb}}

\def\lra{\longrightarrow}


\topmatter

\title  characterization of abelian varieties \endtitle

\author Jungkai A. Chen, Christopher D. Hacon   \endauthor

\address Jungkai Alfred Chen, Department of Mathematics, National Chung
Cheng University, Ming 
Hsiung, Chia Yi, 621, Taiwan \endaddress

\email  jkchen\@math.ccu.edu.tw  \endemail

\address Christopher Derek Hacon, UTAH, Department of Mathematics, Salt
Lake City, UT 84112, USA 
\endaddress

\email   chhacon\@math.utah.edu   \endemail

\footnote""{The first author was partially supported by
National Science Council (NSC-87-2119-M-194-007)}
\abstract
We prove that any smooth complex projective variety $X$
with plurigenera $P_1(X)=P_2(X)=1$ and
irregularity $q(X)=\DI (X)$ is birational to an
abelian variety.
\endabstract

\date March 15, 1999 \enddate

\endtopmatter

\document

\heading Introduction \endheading
Given a smooth projective variety $X$,
one would like to characterize it by its birational invariants.
Classical results include Castelnuovo's criterion for rational surfaces,
which states that a complex surface with  plurigenera $P_2(X)=0$
and irregularity $q(X)=0$ must be rational.
Another well known result 
due to Enriques,
characterizes abelian surfaces as the only surfaces with
$P_1(X)=P_2(X)=1$ and $q(X)=2$.

In an attempt to generalize Enriques' result to
higher dimensional varieties,
Kawamata proved that a variety $X$ with Kodaira dimension
$\kappa(X)=0$
and irregularity $q(X)=\DI (X)$
admits a birational map to an abelian variety \cite{Ka}.
In \cite{Ko1}, Koll\'ar gives an effective characterization by proving that
if
$P_1(X)=P_4(X)=1$ and
$X$ has maximal Albanese dimension,
then $X$ is birational to an abelian variety.
Ein and Lazarsfeld prove a similar result by using generic
vanishing theorems \cite{EL1}.
Koll\'ar subsequently improved his result by
weakening the hypothesis to $P_3(X)=1$ and $q(X)=\DI (X)$.
He conjectured that the result still holds if one assumes that
$P_2(X)=1$ and $X$ is
generically finite over an abelian variety \cite {Ko2}.

In this paper,
we verify Koll\'ar's conjecture.
Our method is based on the approaches of Koll\'ar and of Ein and Lazarsfeld.

{\bf Acknowledgment.} We are in debt to R. Lazarsfeld, L. Ein and J.
Koll\'ar for valuable
conversations during the preparation of this paper.

\heading Conventions and Notations   \endheading

(0.1) Throughout this paper, we work
over the field of complex numbers $\CC$.

(0.2) For $D_1,D_2$ $\QQ$-divisors on a variety $X$, we
write $D_1 \prec D_2$ if $D_2-D_1$ is effective,
and $D_1 \equiv D_2$ if $D_1$ and $D_2$ are numerically equivalent.

(0.3) $|D|$ will denote the linear series associated to the divisor
$D$, and $Bs|D|$ denotes the base locus of $|D|$.

(0.4) For a real number $a$, let $\lfloor a \rfloor$ be the
largest integer $\le a$
and $\lceil a \rceil$ be the smallest integer $\ge a$.
For a  $\QQ$-divisor $D=\sum a_i D_i$, let
$\lfloor D \rfloor =\sum \lfloor a_i \rfloor D_i$ and
$\lceil D \rceil =\sum \lceil a_i \rceil D_i$.

(0.5) Let $\FF$ be a coherent sheaf on $X$,
then $h^i(X,\FF)$ denotes the complex dimension of $H^i(X,\FF)$.
In particular,
the plurigenera $h^0(X,\ox^{ \otimes m})$ are denoted  by $P_m(X)$
and the irregularity $h^0(X,\Omega_X^1)$ is denoted by $q(X)$.

(0.6) Let $A$ be an abelian variety,
for all positive integers $n$,
we will denote by $A_n$ the set of
$n$-division points of $A$.
For all subgroups $T$ of $A$,
$T_0$ will be
the connected component of $T$ containing the origin.
For any line bundle
$L$ there is a canonical homomorphism
$\phi _L :A\longrightarrow \Pic ^0(A)$
defined by $x\lra t_x^*L\otimes L^{-1}$.
Denote by $K(L)$ the kernel of
$\phi _L$, and by $p_L:A\lra A(L):=A/K(L)_0$
the induced quotient map.

\heading 1. Cohomological Support Loci \endheading
We start by recalling  some results about
cohomological support loci that will be needed in this paper.
Most of these results may be found in \cite{EL1} and \cite{EV}.

Let $f:X\to A$ be a morphism from a smooth projective variety $X$ to
an abelian variety $A$.
If $\FF$ is a coherent sheaf on $X$,
then one can define the cohomological support loci by

\proclaim{Definition 1.1}
$$V_i(X,A,\FF):=\{ P \in \Pic ^0(A) | h^i(X, \FF \otimes f^*P) \ne 0
\}.$$
In particular,
if $f={ {\alb }}_X :X \to \Alb (X)$,
then we simply write
$$V_i(X,\FF):=\{ P \in \Pic ^0(X) | h^i(X, \FF \otimes P) \ne 0 \}.$$
\endproclaim

We say that $X$ has maximal Albanese dimension if $\DI {(\alb _X(X))}=\DI
(X)$.
The loci $V_i(X,\ox )$ defined above are very useful
in the study of irregular varieties
and in particular of varieties
of maximal Albanese dimension.
The geometry of these loci is governed by the following:

\proclaim{Theorem 1.2 (\cite{GL1},\cite{GL2})}

a. Any irreducible component of $V_i(X,\ox )$ is a translate of
a sub-torus of $\Pic ^0(X)$
and is of codimension at least
$i-(\DI (X) -\DI (\alb _X(X)))$.

b. Let $P \in S$ be a general point of an irreducible component
$T$ of $V_i(X,\ox)$.
Suppose that $v
\in H^1(X, \OO
_X)\cong T_{P }\Pic ^0(X)$ is not tangent to $T$. Then the sequence
$$H^{n-i-1}(X,\ox \otimes P) \overset {\cup v} \to {\longrightarrow}
H^{n-i}(X,\ox \otimes P)  \overset {\cup v} \to {\longrightarrow }
H^{n-i+1}(X, \ox \otimes P) $$
is exact.
If $v$ is tangent to $T$, then the maps in the above sequence vanish.

c. If $X$ is a variety of maximal Albanese dimension, then
$$\Pic ^0(X)\supset V_0(X,\ox )\supset V_1(X,\ox)
\supset \ ...\ \supset V_n(X,\ox )=\{ \OO _X\}.$$
\endproclaim

\proclaim{Theorem 1.3 (\cite{S})}
Under the above hypothesis,
every irreducible component of $V_i(X,\ox)$ is
a translate of a sub-torus of $\PX$ by a torsion point.
\endproclaim

\remark{Remark 1} The results listed above also hold when we
consider an abelian variety $A$ and a morphism $f: X \to A$ with
$\DI f(X)=\DI (X)$ and we replace $\PX$ by $\Pic ^0(A)$.
\endremark

\proclaim{Proposition 1.4 \cite{EL1, Proposition 2.1}} If
$P_1(X)=P_2(X)=1$, then the origin is an isolated point of
$V_0(X,\omega _X)$.
\endproclaim
\proclaim{Proposition 1.5 \cite{EL1, Proposition 2.2}} If the
origin is an isolated point of $V_0(X,\omega _X)$, then the
Albanese mapping $\alb _X:X\lra \Alb (X)$ is surjective.
\endproclaim

If $X$ has maximal Albanese dimension, then $X$ admits a
generically finite map to an abelian variety $A$. It is easy to see
that $P_1(X) \ge 1$ by pulling back sections of $\Omega^n_A$, where
$n=\dim (X)$. Assume furthermore that $P_2(X)=1$, then necessarily $P_1(X)=1$,
and so the Albanese map is surjective. In particular $q(X)=\dim (X)$.
Therefore one has:

\proclaim{Corollary 1.6} The following are
equivalent.

a. $P_1(X)=P_2(X)=1$, and $q(X)=\dim (X)$.

b. $P_2(X)=1$ and $X$ has maximal Albanese dimension.
\endproclaim

Next we illustrate some examples in which the geometry of $X$ can
be recovered from information on the loci $V_i(X, \ox)$.

\proclaim{Theorem 1.7 (\cite{EL2})}
If $X$ is a variety with
maximal Albanese dimension and $\DI V_0(X,\ox)=0$. Then $X$ is
birational to an abelian variety.
\endproclaim

The idea of the proof is as follows. There exists a map $f:X\lra
A$ to an abelian variety $A$, such that $X$ is generically finite
over its image. We may assume that the map $f:X\lra A$ is
surjective, as otherwise by Proposition 1.5, there exists a
positive dimensional component of $V_0(X,\omega _X)$. Consider the
sheaf $f_*\omega _X$. By Theorem 1.2, the groups $H^i(A, f_*\omega
_X\ot P)$ do not vanish only for finitely many $P\in \Pic ^0(A)$.
(In fact it is not hard to see \cite{EL2} that $h^i(f_*\omega _X\ot P)=0$
for all $P\not \hskip-.12cm{=}\OO _X$.)
By \cite{M, Example 2.9 and \S 3}, $c_1(f_*\omega _X)=0$, and so
the map $X\lra A$ is \'etale in codimension one and hence \'etale.

\proclaim{Corollary 1.8 (\cite{EL2})} Let $X$ be a variety of
maximal Albanese dimension, and of positive Kodaira dimension ($\kappa(X)>
0$). Then there exists a torsion line bundle $Q\in \Pic ^0(X)$ and
a positive dimensional subgroup $T$ of $\Pic ^0(X)$, such that
$h^0(2K_X\ot 2Q\ot P)\geq 2$ for all $P\in T$.
In particular, if $V_0(X,\ox)$ contains a positive dimensional
component through the origin $\OO _X$ (or through any 2-torsion
point $P\in \Pic ^0(A)_2$), then $P_2(X) \ge 2$.
\endproclaim

\demo{Proof} Let $S_Q$ be any positive dimensional component of
$V^0(X,\omega _X)$ containing
$Q \in \Pic ^0(X)$ and let $S=S_Q-Q$. Since $S$ is a subgroup of
$\Pic ^0(X)$, one sees that $-S=S$. Consider
the natural map 
$$H^0(X,\ox \ot Q \ot P_y) \otimes H^0(X, \ox \ot
Q \ot P_y^*) \to H^0(X, \ox^{ \otimes 2} \ot Q^{ \otimes 2})$$ 
and let $P_y$
vary in $S$. Both $H^0(X,\ox \ot Q \ot P_y)$ and $H^0(X,\ox \ot Q
\ot P_y^*)$ are non-zero. The natural map above gives rise to a
map of linear series $$ |K_X +Q+P_y| \times |K_X+Q -P_y| \to
|2K_X+2Q |.$$ Since there are infinitely many $P_y \subset S$,
$|2K_X+2Q|$ cannot consist of only one element. Therefore
$h^0(X,\ox^{\ot 2}\ot Q^{ \otimes 2})\ge 2$. \hfill \qed
\enddemo

\proclaim{Corollary 1.9 \cite{EL1}} Let $D$ be a reduced
irreducible divisor of an abelian variety $A$. Then there exists a
positive dimensional sub-group $T_D\subset \Pic ^0(X)$ such that
for any desingularization $\nu:\tilde {D}\lra D$, $h^0(\omega
_{\tilde {D}} \ot \nu ^* P)>0$ for all $P\in T_D$.
\endproclaim
We end this section by recalling a few results that will 
frequently be used in what follows.
\proclaim{Theorem 1.10 \cite  {Ko1, EV2}} 
Let $f:X\lra Y$ be a surjective morphism from a smooth 
projective variety $X$ to a normal variety $Y$. Let $L$ be a line bundle on $X$
such that $L\equiv f^* M +\Delta$, where $M$ is a $\QQ$-divisor on $Y$
and $(X,\Delta )$ is klt. Then

a. $R^jf_*(K_X\ot L) $ is torsion free for $j\geq 0$.

b. Assume furthermore that $M$ is nef and big. Then $H^i(Y,R^jf_*(K_X\ot L))
=0$ for $i>0$, $j\geq 0$.
\endproclaim

\proclaim{Theorem 1.11 \cite  {Ko2}} 
Let $f:X\lra Y$ be a surjective morphism from a smooth 
projective variety $X$ to a smooth variety $Y$, birational 
to an abelian variety.
Let $L$ be a nef and big $\QQ$-divisor on $Y$, $\Delta$ be a $\QQ$-divisor
with normal crossings such that $\lfloor \Delta \rfloor =0$. Let $U\subset 
Y$ be a dense open set, and $N$ a Cartier divisor on $X$, $N\equiv K_X+\Delta
+f^*L$. Assume that $N|g^{-1}U$ is linearly equivalent to an effective
divisor. Then $h^0(X,N)\not \hskip-.12cm {=}0$.
\endproclaim 

\proclaim{Theorem 1.12 \cite {Ko2}} Let $f:X\lra Y$ be a
surjective morphism from a smooth projective variety $X$ to a
normal variety $Y$. Let $L$ be a line bundle on $X$ such that
$L\equiv f^* M +\Delta$, where $M$ is a $\QQ$-divisor on $Y$ and
$(X,\Delta )$ is klt. Let $D$ be a effective divisor on $X$
such that $f(D) \ne Y$. Then

$$H^j(X,K_X \ot L) \to H^j(X, K_X \ot L(D)) \text{ is
injective.}$$
\endproclaim

\remark{Remark} Suppose $X$ admits a generically finite map $a: X
\to A$ to an abelian variety $A$. And let $f: X \to Y$ be its
Iitaka fibration. Let $D$ be an non-exceptioanl irreducible
component of the ramification divisor with respect to the
map $a: X \to A$. Then
$P_1(D) \ge 1$, and $f(D) \ne Y$ \cite{Ko2, 17.6.1}.
\endremark

\heading 2. Birational invariants \endheading In this section we
study  birational invariants of varieties with maximal Albanese
dimension. 
Let $f:X \to Y$ be a birational model of the Iitaka fibration of $X$.
We may assume that $X$ and $Y$ are
smooth projective
varieties.
Let $K$ be the image of the generic fiber of $f$ via the Albanese map
$\alb _X:X\lra \Alb (X)$.
K$ $ is a sub-group of $\Alb (X)$, and
there is a generically finite map $Y \longrightarrow \Alb (X)/K$.
Thus $Y$ is also of maximal Albanese dimension \cite{Ko2, 17.5.1}.

\proclaim{Lemma 2.1} Suppose that $X$ has maximal Albanese
dimension, and positive Kodaira dimension $\kappa(X) >0$. Fix $Q$ a
torsion element of $ \Pic ^0(X)$. Then $h^0(X,2K_X \ot Q \ot
f^*P)$ is constant for all torsion $P \in \Pic ^0(Y)$.
\endproclaim

\demo{Proof} Fix a nef and big divisor $A$ on $Y$. The linear
series $|mK_X-f^*A|$ is nonempty for all $m$ sufficiently large
\cite{Ko2, Lemma 17.6.1}. Fix an appropriate divisor
$B_2 \in |m_2K_X-f^*A|$. For an arbitrary integer $m_1$, one may
consider the map 
$$|m_1K_{X}|\lra |(m_1+m_2)K_{X}-{f}^*A|,$$
induced by multiplication by $B_2\in |m_2K_{X}-{f}^*A|$. 
Replacing $X$ by an appropriate birational model, we may assume that for
all sufficiently big and divisible integers $m_1$,
$|m_1K_X|=|M_{m_1}|+F_{m_1}$ and $|m_2K_X-f^*A|=|M_{m_1,A}|+F_{m_1,A}$,
where $M_{m_1}$ and $M_{m_1,A}$ are free and $F_{m_1}+F_{m_1,A}$ is a divisor
with normal crossings. In particular we may assume that
the divisors $F_{m_1}=Bs |m_1K_{X}|$ and $F_{m_1+m_2,A}=
Bs|(m_1+m_2)K_{X}-{f}^*A|$
have normal crossing singularities.
It
follows that 
$$ \lfloor \frac {Bs|(m_1+m_2)K_{X}-{f}^*A|}
{m_1+m_2} \rfloor \prec \lfloor \frac {Bs |m_1K_{X}|+B_2}{m_1+m_2}
\rfloor .$$
We will show next that for all sufficiently big and divisible integers $m_1$,
$$\lfloor \frac {Bs
|m_1K_{X}|+B_2}{m_1+m_2} \rfloor \prec \frac {2
Bs|m_1K_X|}{m_1}.$$

We proceed by comparing the multiplicities of each component. One may
assume that for all sufficiently big and divisible integers $m_1$,
the support of $Bs|m_1K_X|$ and of $B_2$, are contained in a fixed divisor $D$
with irreducible components $D_i$, $i=1 \ldots r$. (Fix for example
$m_1=2^jm_0$, $j\geq 0$.) We 
may therefore write
$Bs|m_1K_X|=\sum_{i=1}^r b_{m_1,i} D_i$, and $B_2=\sum_{i=1}^{r}
\delta_{i} D_i$.
If $b_{m_1,i}$ is bounded, then $\lfloor \frac
{b_{m_1,i}+\delta_i}{m_1+m_2} \rfloor=0$, and hence is bounded by the 
positive rational numbers $\frac{
2b_{m_1,i}}{m_1}$ for all $m_1$ sufficiently big and divisible.
Conversely, assume that $b_{m_1,i}$ is not bounded. Then for all 
$m_1$ sufficiently big and divisible
$$ \frac
{b_{m_1,i}+\delta_i}{m_1+m_2} < \frac {2b_{m_1,i}}{m_1+m_2} <
\frac {2b_{m_1,i}}{m_1}.$$ 
Therefore the desired inequality holds.

Since $Q$ and $P$ are torsion, we may assume that $m_1=2m_1'$ is
chosen so that $m_1'(f^*P)=m_1'Q=\OO _{X}$, and therefore

$$\frac {2Bs|m_1K_{X}|} {m_1}=\frac {Bs|m_1'(2K_{X}\ot Q \ot f^*P
)|} {m_1'} \prec Bs |2K_{X}\ot Q \ot f^*P | .$$ Therefore one can
choose a divisor $B \in |(m_1+m_2)K_X-f^*A|$ such that

$$\lfloor \frac {B}{m_1+m_2} \rfloor =\lfloor \frac {Bs
|(m_1+m_2)K_X-f^*A|}{m_1+m_2} \rfloor \prec Bs |2K_{X}\ot Q \ot
f^*P | ,$$
and such that
$$\lfloor \frac {B}{m_1+m_2} \rfloor \prec Bs |2K_{X}\ot Q|.$$
Now let $m=m_1+m_2$, and define
$$M:= K_X(- \lfloor \frac{B}{m} \rfloor) \equiv f^*
(\frac{A}{m}) + \{ \frac{B}{m} \}.$$
One sees that $h^0(X,K_X+M+Q+f^*P)=h^0(X,2K_X+Q+f^*P)$, and
$h^0(X,K_X+M+Q)=h^0(X,2K_X+Q)$.
Since the support of $B$ has only normal crossing singularities, 
$\{ \frac{B}{m} \}$ is
Kawamata log terminal.
By Theorem 1.10,
it follows that $H^i(Y,R^j{f}_*(K_{X}\ot M \ot T)\ot
P)=0$ for all $i>0$, $j\geq 0$ and all $T \in \Pic ^0(X)$. In
particular
$$h^0(Y,{f}_*(K_{X}\ot M \ot Q)\ot P)=\chi(Y,{f}_*(K_{X}\ot M \ot
Q)\ot P),$$
which is deformation invariant. So, the quantity
$h^0(X,K_X \ot M \ot Q \ot f^*P)$ is constant for all $P \in
\Pic ^0(Y)$.
Therefore, for any torsion element $P$ of $\Pic ^0(Y)$,
$$h^0(X,2K_X \ot Q \ot f^*P)=h^0(X,K_X \ot M \ot Q \ot
f^*P)=$$
$$h^0(X,K_X \ot M \ot Q)=h^0(X,2K_X  \ot Q).$$
\hfill \qed
\enddemo

An immediate consequence of Lemma 2.1 is the following

\proclaim{Theorem 2.2} Let $X$ be a smooth projective variety with
maximal Albanese dimension. If $X$ is of general type, then
$P_2(X)= h^0(X, \ox^{\ot 2}) \ge 2$.
\endproclaim

\demo{Proof} Assume that $P_2(X)=1$, then by Lemma 1.6, the map
$X\lra \Alb (X)$ is surjective and generically finite. Since $X$
is not birational to an abelian variety, by Corollary 1.9, there
exists a torsion line bundle $Q \in \Pic ^0(X)$, such that
$h^0(X,2K_X+2Q) \ge 2$. By Lemma 2.1, $h^0(X,2K_X) \ge 2$, and
this is the required contradiction. \qed
\enddemo

\proclaim{Lemma 2.3} Keeping the notation as in Lemma 2.1. Assume
furthermore that the Albanese map $g:X\lra \Alb (X)$ is surjective and
generically finite with ramification divisor $R$. Let $R_0$ be any
reduced irreducible smooth divisor contained in $R$ with
multiplicity $r_0\geq 2$, then $h^0(X,2K_X-(r_0-1) R_0 \ot
f^*P)\geq 1$ for all $P \in \Pic ^0(Y)$. Moreover, if the divisor
$R_0$ is not exceptional with respect to the map $g:X\lra \Alb
(X)$, then $P_2(X)\geq 2$.
\endproclaim

\demo{Proof} The proof is a technical improvement on a result of Koll\'ar.
We closely follow the methods and the notation of \cite {Ko2, \S 17}.
There is no harm in assuming that $R_0$ be a smooth divisor as this
can always be arranged replacing $X$ by an appropriate birational model.

Consider the divisor
$$K'=K_X-(r_0-1)R_0,$$
we have the following inclusion of effective divisors
$$K_X\prec r_0K'\prec r_0K_X.$$
It follows by \cite{Ko2, 17.6.1}
that for any fixed ample divisor $A$ on $Y$, the linear series
$|mK'-f^*A|$ is not empty, for all sufficiently big and divisible
$m$. Let $B'$ be a general member of the above linear series, then
define
$$M':=K'(-\lfloor (1/m)B'\rfloor )\equiv f^*(\frac {A}{m})+\{ \frac
{B'}{m}\} .$$
We may further assume that $R_0$ is not a component of $B'$
as by \cite{Ko2}, one can arrange that $R_0$ is not contained
in the base locus of $|mK_X-f^*A|$. However,
$$mr_0K'-f^*A=m(r_0K'-K_X)+mK_X-f^*A,$$
and $r_0K'-K_X$ is an effective divisor not containing $R_0$.

By Theorem 1.10, $H^i(Y,f_*(K_X\ot M')\ot P)=0$
for all $i>0$ and all $P\in \Pic ^0(Y)$. Hence the quantity
$h^0(f_*(K_X\ot M')\ot P)$ is positive and independent of $P\in \Pic ^0(Y)$. 
Therefore $h^0(X,2K_X-(r_0-1) R_0 \ot f^*P)\geq 1$ for
all $P \in \Pic ^0(Y)$.

Consider the exact sequence \cite {Ko2, 17.7},
$$0\lra K_X\ot M'\lra K_X\ot M'(R_0)\lra K_{R_0}\ot M'\lra 0.$$
By Theorem 1.12, the map $H^1(X,K_X\ot M')\lra H^1(X,K_X\ot M'(R_0))$
is injective.
Since, $M'$ and $M'|R_0\equiv f^*A|R_0+\{ B/m|R_0\}$ satisfy the assumptions
of
Theorem 1.11, it follows that the quantities 
$h^0(K_X\ot M')$, $h^0(K_{R_0}\ot M')$ are both positive, 
and hence $h^0(K_X+M'(R_0))\geq 2$.
The assertion now follows from the map of linear series $|K_X\ot M'(R_0)|
\lra |2K_X|$ induced by multiplication by the effective divisor
$(r_0-2)R_0+\lfloor (1/m)B'\rfloor$.
\hfill \qed
\enddemo

\proclaim{Theorem 2.4} Let $X$ be a smooth projective variety with
maximal Albanese dimension and positive Kodaira dimension $\kappa(X) >0$.
Let $f:X\to Y$ be
the Iitaka fibration. If $Y$ is birational to an abelian variety,
then $P_2(X) \ge 2$.
\endproclaim

\demo{Proof}  Let $q:Y \to S$ be a birational morphism to an
abelian variety $S$. In what follows, we assume that $h^0(X,2K_{X}) =1$
and derive a contradiction.

{\bf Step 1.} Since $X$ is  not birational to an abelian variety,
by Corollary 1.8,
we may choose
a torsion line bundle $Q$ such that
$h^0(X,2K_X+2Q) \ge 2$. By Lemma 2.1, $h^0(X,2K_{X}+2Q+P) \ge
2$ for all torsion $P \in \Pic ^0(Y)$. It follows by the upper
semicontinuity
of the dimension of cohomology groups, that $h^0(X,2K_{X}+2Q+P) \ge  2$
for all $P \in \Pic ^0(Y)$.

{\bf Step 2.} Again by Lemma 2.1, since the dimension of
cohomology groups is upper semicontinuous, $h^0(X,2K_{X}+P) \ge  1$
for all $P \in \Pic ^0(Y)$. As $h^0(2K_X)=1$, equality holds for all $P$ in
an
appropriate neighborhood $U\subset \Pic ^0(Y)$ of the origin $\OO _Y$.

{\bf Step 3.} There exists an effective divisor
$G$ which is contained in the base locus of $|2K_X+f^*P|$ for
all $P\in U$.

To this end, let $M$ be a divisor chosen as in the proof of Lemma
2.1 (in particular $\lfloor \frac{B}{m} \rfloor \prec Bs|2K_X|$).
We wish to compute the cohomology groups of the direct image sheaf
$q_*f_*(K_{X}+M) $.
From the Leray
spectral sequence associated to the morphism $q:Y \lra S$, one sees that
for all $P\in \Pic ^0(Y)$, 
$$h^0(X,K_{X}\ot M\ot f^*P)=h^0(S,q_*f_*(K_{X}+M) \ot P )=1,$$
and there is an injection
$$H^i(S,q_*{f}_*(K_{X}\ot M)\ot P )\hookrightarrow H^i(Y,{f}_*(K_{X}\ot M)
\ot P).$$
It follows that for all $i>0$,
both of the above cohomology groups vanish.
So the sheaf
$q_*{f}_*(K_{X}\ot M)$ is coherent, generically of rank one
which is a "cohomological principal polarization". In other words,
for all $P\in \Pic ^0(S)$, the
cohomology groups of $q_*{f}_*(K_{X}\ot M)\ot P$
are of the same dimension as
the cohomology groups of a principal polarization of $S$. By
\cite{H, Proposition 2.2},
there exists an appropriate theta divisor $\Theta$ on $S$ and
an isomorphism of sheaves $q_*{f}_*(K_{X}\ot M)
 = \OO _S(\Theta )$.
Pulling back with respect to $q$ and $f$, for all $P\in \Pic ^0
(Y)$, we get a map of linear series
$${f}^*q^*|\OO _S(\Theta )\ot P|@>{\times G}>> |2K_{X} \ot
f^*P|$$
which factors through $|K_{X}\ot M\ot f^*P|$.
Since for
$P\in U$, the groups
$H^0(X,2K_{X} \ot f^*P)$ are one dimensional, the above map of linear
series is an isomorphism, and so $G$ is contained in $Bs|2K_{X}
\ot f^*P|$ for all $P\in U\subset \Pic ^0(Y)$.

{\bf Step 4.} We wish to show that the divisor $G$ is non-zero and
not exceptional with respect to the map $X\lra S$.

Let $F$ be a maximal effective divisor such that for all $P\in \Pic ^0(Y)$,
$2K_X\ot 2Q\ot f^* P$ is linearly equivalent to $F+V_P$,
and $h^0(V_P)\geq 2$.
%
%
We may assume that
$D_P,\ D'_P$ are two distinct general members of $|V_P|$, which do
not contain any components of the ramification divisor of $X\lra
Y$. Therefore there exist divisors $\Delta _P,\Delta '_P$ on $Y$
such that $D_P={f}^*(\Delta _P)$ and $D'_P={f}^*(\Delta '_P)$. The map $Y
\lra S$ is birational, so there exist effective divisors $\bar {\Delta
}_P$ and $\bar {\Delta }'_P$ on $S$ and arbitrary $q$-exceptional
divisors $E_P,\ E'_P$ such that

$$\Delta _P=q^*(\bar {\Delta }_P)+E_P,\ \ \ \Delta '_P=q^*( \bar
{\Delta }'_P)+E'_P.$$
If $h^0(S,\OO _S(\bar {\Delta }_P))=1$, then $\bar {\Delta }_P=
\bar {\Delta }'_P$ and similarly since $Y \lra S$ is birational,
$\Delta _P=\Delta '_P$, which is a contradiction. So $\bar {\Delta
}_P$ is a divisor on $S$ with at least two sections. We may also
assume that $\bar {\Delta }_P$ is an ample divisor on $S$, since
otherwise there would exists a non trivial map of abelian
varieties $S\lra \Sigma$ such that $\bar {\Delta }_P$ is the pull
back of a divisor on $\Sigma$. However then for general $P_1,$
$P_2\in \Pic ^0 (S)$, $V_{P_1}-V_{P_2}$ is rationally equivalent to
$P_1-P_2$ and is an element of $\Pic ^0 (\Sigma )$. The only way
this can occur is if $S=\Sigma$ as required.

Hence, for $H'$ the pull back of an ample divisor on $A=\Alb (X)$, and $H$
the pull back of an ample divisor on $S$, we have the numerical inequality
$$2K_{X}\cdot H^{s-1}\cdot (H')^{n-s}> (f^*q^*\Theta) \cdot
H^{s-1} \cdot (H')^{n-s}.$$
From the injection of linear series
$${f}^*q^*|\OO _S(\Theta )\ot
P|@>{ \times G}>> |2K_{X} \ot f^*P|,$$ one sees that the effective
divisor $G$ is not exceptional with respect to the map
$X\lra S$.

{\bf Step 5.} Consider now $G_0$ an irreducible reduced component
of $G$ which is not exceptional with respect to the map $X\lra S$.
Since we assumed that
$h^0(X,2K_X)=1$, by Proposition 1.4,
$h^0(X,K_X \ot P)$ $=0$ for $P\in (\Pic ^0(Y)- \OO _Y)$ in
an appropriate
neighborhood $U'$ of the origin $\OO _S$. We may assume that $U=U'$.
From the exact sequence of sheaves
$$0\lra K_X \ot P\lra K_X+G_0\ot P\lra K_{G_0} \ot P\lra 0,$$
we see that for $P\in U^0=U-\{\OO _Y\}$, there is an isomorphism
$$H^0(X,K_X+G_0\ot P)\cong H^0(G_0,K_{G_0} \ot P).$$
By Corollary 1.9,
there exists a positive dimensional subgroup $T_{G_0}\subset
\Pic ^0(X)$ such that for all $P\in T_{G_0}$, $h^0(K_{G_0} \ot P)>0$.
Since $G_0\subset 2K_X$, $G_0$ is vertical with respect to the
Iitaka fibration
$X\lra Y$. We may therefore assume that $T_{G_0}$ is in fact contained in
$\Pic ^0(Y)$.
It follows that for $P\in T_{G_0}\cap U^0$, $G_0$ is not contained in the
base locus of the linear series $|K_X+G_0 \ot P|$.
Let
$$|K_X+G_0 \ot P| @>{K_X-G_0}>> | 2K_X \ot P|$$
be the map of linear series induced by multiplication by the
effective divisor $K_X -G_0$.
Since $G_0$ is in the base locus of $|2K_X\ot P|$, $G_0$
must be an irreducible component of $K_X -G_0$.
In other words the multiplicity $g_0$ of $G_0$ in $R=K_X$ is at least $2$.
By Lemma 2.3, this is the required contradiction.
\qed

\enddemo

\heading 3. Proof of the Main Theorem \endheading

\proclaim{Lemma 3.1}
Let $X$ and $Y$ be smooth varieties and
let $A$ and $C$ be abelian varieties.
Suppose that there is a commutative diagram
$$
\CD
X &@>{g}>>&A \\
@V{f}VV && @VV{\pi}V \\
Y & @>{q}>>& C,\\
\endCD
$$
such that both $g$ and $q$ are generically finite and surjective.
Then $P_m(Y) \le P_m(X)$, for all $m \ge 1$.
\endproclaim

\demo{Proof}
Let
$\{ \bar{z}_1,\ldots,\bar{z}_k\}$ be coordinates on $C$
such that
$\{z_i= \bar{z}_i\circ \pi \ |\ 1\leqq i\leqq k\ \}$
may be extended to coordinates
$\{ z_1,\ldots,z_n  \}$ on $A$,
where $n$ and $k$ are the dimensions of $A$ and of $C$ respectively.
Let $\Lambda=\WED {g^*dz} {k+1}
n \in H^0(X, \Omega ^{n-k}_{X})$,
we claim that
$$ (f^* \omega_Y)^{\otimes m} @>{\wedge \Lambda ^{\otimes m}}>>
\omega_X^{\otimes m}$$
induces the desired injective map
$$\Lambda _m: H^0(Y,\omega_Y^{\otimes m}) @>>> H^0(X,\omega_X^{\otimes
m}).$$
Since $q:Y \to C$ is generically finite and surjective,
for a general point $y \in Y$,
there are isomorphic Euclidean open neighborhoods
$\UU_y$ and $\VV_{q(y)}$ of $y$ and $q(y)$ respectively.
For a non-zero section
$\eta \in H^0(Y,\omega_Y^{\otimes m})$,
we may assume that $\eta(y) \ne 0$.
On $\UU_y$, $\eta$ is represented by
$\zeta\cdot (\WED {q^*d\bar{z}} 1 k)^m$, with $\zeta (y) \ne0$.
On $f^{-1}(\UU_y)$
$$ \Lambda _m( \eta) |_{f^{-1}(\UU_y)} =
f ^*\zeta \cdot (g^*dz_1 \wedge  \ldots \wedge g^*dz_n)^m\not
\equiv 0.$$
Therefore
(by choosing harmonic representatives, see \cite{GL2, Lemma 3.1}),
$\Lambda _m: H^0(Y,\omega_Y^{\otimes m}) @>>> H^0(X,\omega_X^{\otimes
m})$
is injective and so $P_m(Y) \le P_m(X)$.
\hfill \qed
\enddemo

\proclaim{Theorem 3.2}
Let $X$ be a smooth projective variety with
$P_1(X)=P_2(X)=1$ and $q(X)=\DI (X)$.
Then $X$ is birational to an abelian variety.
\endproclaim

\demo{Proof}
Since $P_1(X)=P_2(X)$ and $q(X)=\DI (X)$, by Proposition 1.5,
it follows that the Albanese morphism is surjective, and hence
a generically finite map.

If the Kodaira dimension of $X$ is zero, then
by Kawamata's theorem \cite {Ka}, $X$ is birational to an abelian variety.

Suppose that $\kappa(X) >0$.
By Theorem 2.2, the variety $X$ cannot be of general type.
Therefore $X$ admits a nontrivial Iitaka fibration.
Let $f: X \to Y$ be a birational model of the Iitaka fibration.
The generic fiber $F$ has
$\kappa(F)=0$ and $\DI ({\alb _X(F)}) =\DI (F)$.
Thus $\alb _X(F)$ is an abelian sub-variety which we denote by $K$.
We obtain a diagram
$$
\CD
 X & @>{ \alb _X}>> & \Alb (X) \\
@V{f}VV & & @VV{\pi}V \\
Y & @>{g}>> & \Alb (X)/K\\
\endCD
$$
where the morphism $g: Y @>>> \Alb (X)/K$ is generically finite
(cf [Ko2, 17.5.1]).
By Lemma 3.1, we have $P_1(Y)=P_2(Y)=1$.
By iterating this procedure,
we eventually produce a variety $Y'$
either of maximal Kodaira dimension $\kappa (Y') = \DI (Y')$,
or with Kodaira dimension $0$ (and hence birational to an abelian variety).
By applying Theorem 2.2 or Theorem 2.4,
respectively, in conjunction with Lemma 3.1,
we conclude that $P_2(X) \ge 2$.
\hfill
\qed
\enddemo


\enddocument
\end